\documentclass[12pt]{amsart}
\usepackage{latexsym,amssymb,amscd}

\theoremstyle{plain}
   \newtheorem{theorem}{Theorem}
   \newtheorem{proposition}[theorem]{Proposition}
   \newtheorem{lemma}[theorem]{Lemma}

\theoremstyle{definition}

\newenvironment{proofof}[1]{\begin{trivlist}\item {\it
        Proof of {#1}.\,}}{\mbox{}~\hfill$\Box$\end{trivlist}}

\newenvironment{numbered}%
        {\begin{list}
                {\noindent\makebox[0mm][r]{\arabic{enumi}.}}
                {\leftmargin=5.5ex \usecounter{enumi}\topsep=1.5mm}
        }
        {\end{list}}

\newenvironment{bullets}%
        {\begin{list}
                {\noindent\makebox[0mm][r]{$\bullet$}}
                {\leftmargin=5.5ex \usecounter{enumi} \topsep=1.5mm}
        }
        {\end{list}}

\numberwithin{equation}{section}

\renewcommand{\>}{\rangle}

\newcommand\AAA{{\mathcal A}}
\newcommand\RR{{\mathbb R}}
\newcommand\ZZ{{\mathbb Z}}
\newcommand\QQ{{\mathbb Q}}

\newcommand\PL{\textit{PL}}
\newcommand\FF{{\mathcal F}}
\newcommand\GG{{\mathcal G}}

\newcommand\kk{{\Bbbk}}
\newcommand\mm{{\mathfrak m}}
\newcommand\Hilb{{\mathrm{Hilb}}}
\newcommand\Ext{{\mathrm{Ext}}}
\newcommand\Hom{{\mathrm{Hom}}}
\newcommand\link{{\mathrm{link}}}

\newcommand\rank{{\mathrm{rank}}}
\newcommand\minus{\smallsetminus}
\newcommand\nullity{{\mathrm{nullity}}}

\begin{document}

\title{Reciprocal domains and Cohen--Macaulay $d$-complexes in $\RR^d$}

\author{Ezra Miller}
\author{Victor Reiner}
\email{(ezra,reiner)@math.umn.edu}
\address{School of Mathematics\\
University of Minnesota\\
Minneapolis, MN 55455, USA}

\dedicatory{Dedicated to Richard P. Stanley on the occasion of his 60th birthday}

\thanks{EM and VR supported by NSF grants DMS-0304789 and DMS-0245379
respectively.}


\date{8 April 2004}


\begin{abstract}
We extend a reciprocity theorem of Stanley about enumeration of
integer points in polyhedral cones when one exchanges strict and weak
inequalities.  The proof highlights the roles played by
Cohen--Macaulayness and canonical modules.  The extension raises the
issue of whether a Cohen--Macaulay complex of dimension~$d$ embedded
piecewise-linearly in~$\RR^d$ is necessarily a $d$-ball.  This is
observed to be true for $d \leq 3$, but false for $d=4$.
\end{abstract}

\maketitle

\section{Main results} \label{main}

This note begins by dealing with the relation between enumerators of
certain sets of integer points in polyhedral cones, when one exchanges
the roles of strict versus weak inequalities
(Theorem~\ref{t:reciprocal}).  The interaction of this relation
with the Cohen--Macaulay condition then leads us to study
piecewise-linear Cohen--Macaulay polyhedral complexes of dimension~$d$
in Euclidean space~$\RR^d$ (Theorem~\ref{CM-d-complexes-in-d-space}).

We start by reviewing a result of Stanley on Ehrhart's notion of
reciprocal domains within the boundary of a convex polytope.  Good
references for much of this material are
\cite[Chapter~6]{BrunsHerzog}, \cite[Chapter~1]{Stanley-greenbook},
and \cite[Part~II]{MillerSturmfels}.

Let $Q \subset \ZZ^d$ be a saturated affine semigroup, that is, the
set of integer points in a convex rational polyhedral cone
$C=\RR_{\geq 0}Q$.  Assume that the cone~$C$ is of full dimension~$d$,
and pointed at the origin.  Denote by~$\FF$ the facets (subcones of
codimension~$1$) of $C$.  For each facet $F \in \FF$, let $\ell_F(x)
\geq 0$ be the associated facet inequality, so that the semigroup
\begin{eqnarray*}
  Q &=& \{ x \in \ZZ^d \mid \ell_F(x) \geq 0 \hbox{ for all facets }F
  \in \FF\}
\end{eqnarray*}
is the intersection of the corresponding closed positive half-spaces.

Fix a nonempty proper subset~$\GG$ of the facets~$\FF$, and
let~$\Delta$ and~$\Delta'$, respectively, denote the pure
$(d-1)$-dimensional subcomplexes of the boundary complex of~$C$
generated by the facets in~$\GG$ and $\FF\minus\GG$, respectively.
Ehrhart called the sets $C\minus\Delta$ and $C\minus\Delta'$ {\em
reciprocal domains}\/ within the boundary complex of~$C$.  Examples of
reciprocal domains arise when $\Delta$ is {\em linearly separated from
$\Delta'$}, meaning that some point $p \in \RR^d$ satisfies
$$
  \begin{aligned}
  \ell_F(p) > 0 &\text{ for }F \in \GG \\
  \hbox{and } \ell_F(p) < 0 &\text{ for }F \in \FF\minus\GG.
  \end{aligned}
$$
Define the {\em lattice point enumerator}\/ to be the power series
\begin{eqnarray*}
  F_{C\minus\Delta}(x) &:=& \sum_{a \in (C \minus\Delta)\cap\ZZ^d}x^a
\end{eqnarray*}
in the variables $x=(x_1,\ldots,x_n)$.  This series lies in the
completion~$\ZZ[[Q]]$ of the integral semigroup ring~$\ZZ[Q]$ at the
maximal ideal $\mm=\<x^a \mid 0 \neq a \in Q\>$ generated by the set of nonunit
monomials.  General facts about Hilbert series of finitely
generated modules over semigroup rings imply that
$F_{C\minus\Delta}(x)$ can be expressed in the complete ring~$\ZZ[[Q]]$ as a rational
function whose denominator is a product of terms having the form
$1-x^a$ \cite[Chapter~8]{MillerSturmfels}.

A result of Stanley \cite[Proposition 8.3]{Stanley-reciprocity} says
that when $\Delta,\Delta'$ are linearly separated,
\begin{eqnarray*}
  F_{C\minus\Delta'}(x^{-1}) &=& (-1)^d \, F_{C\minus\Delta}(x)
\end{eqnarray*}
as rational functions in $\QQ(x_1,\ldots,x_d)$.  Our main result
weakens the geometric `linearly separated' hypothesis on~$\Delta$ to
one that is topological and ring-theoretic.

Let $\kk$ be a field, and denote by $\kk[Q] = \bigoplus_{a \in Q} \kk
\cdot x^a$ the $\ZZ^d$-graded affine semigroup ring corresponding
to~$Q$.  For each subcomplex $\Delta$ of~$C$, this ring contains a
radical, $\ZZ^d$-graded ideal~$I_\Delta$ consisting of the $\kk$-span
of monomials~$x^a$ for $a \in C\minus\Delta$.  The {\em face ring}\/
of~$\Delta$ is defined to be the quotient
$\kk[\Delta]:=\kk[Q]/I_\Delta$.

A~polyhedral subcomplex $\Delta \subseteq C$ is {\em Cohen--Macaulay
over $\kk$}\/ if $\kk[\Delta]$ is a Cohen--Macaulay ring.  This turns
out to be a topological condition, as we now explain.  Fix a
$(d-1)$-dimensional {\em cross-sectional polytope}\/~$\overline{C}$ of
the cone $C$, and let $\overline{\Delta}:=\overline{C} \cap \Delta$, a
pure $(d-2)$-dimensional subcomplex of the boundary complex of
$\overline{C}$.  It is known \cite{Yanagawa}
that $\Delta$ is Cohen--Macaulay if and
only if the geometric realization $|\overline{\Delta}|$ is {\em
topologically Cohen--Macaulay (over\/~$\kk$)}, meaning that its
(reduced) homology $\tilde{H}_i(|\overline{\Delta}|;\kk)$ and its
local homology groups
$\tilde{H}_i(|\overline{\Delta}|,|\overline{\Delta}|\minus p;\kk)$
vanish for $i < d-2$.

The Cohen--Macaulay condition is known to hold whenever
$|\overline{\Delta}|$ is a $(d-2)$-ball, but this sufficient condition
is not in general necessary; see
Theorem~\ref{CM-d-complexes-in-d-space} below.  Nevertheless, when
$\Delta$ is linearly separated from $\Delta'$, the topological space
$|\overline{\Delta}|$ is such a ball, because its facets are shelled
as an initial segment of a {\em (Bruggesser--Mani) line-shelling}\/
\cite[Example 4.17]{OMbook} of the boundary complex of the cone $C$.

\begin{theorem} \label{t:reciprocal}
Let $\Delta$ be a dimension~$d-1$ subcomplex of a pointed rational
polyhedral cone~$C \subseteq \RR^d$.  If $\Delta$ is Cohen--Macaulay
over some field\/~$\kk$, then as rational functions, the lattice point
enumerators of the reciprocal domains \mbox{$C\minus\Delta$ and
$C\minus\Delta'$ satisfy}
\begin{eqnarray*}
  F_{C\minus\Delta'}(x^{-1}) &=& (-1)^d \, F_{C\minus\Delta}(x).
\end{eqnarray*}
\end{theorem}

Theorem~\ref{t:reciprocal} raises the issue of whether a
$d$-dimensional Cohen--Macaulay proper subcomplex of the boundary of a
$(d+1)$-polytope must always be a $d$-ball, a question that arises in
other contexts within combinatorial topology (such as
\cite{BilleraRose}).  Although the following result is surely known to
some experts, its (statement and) proof seems sufficiently difficult
to find written down that we include the details in
Section~\ref{d-space}.

\begin{theorem} \label{CM-d-complexes-in-d-space}
\begin{numbered}
\item
Let $K$ be a $d$-dimensional proper subcomplex of the boundary of
$(d+1)$-polytope.  If $d \leq 3$ and $K$ is Cohen--Macaulay over some
field\/~$\kk$, then the topological space $|K|$ is homeomorphic to a
$d$-ball.

\item
There exists a proper subcomplex of dimension\/~$4$ in the boundary of
a $5$-polytope that is Cohen--Macaulay over every field but not
homeomorphic to a $4$-ball.
\end{numbered}
\end{theorem}

\section{Reciprocal domains via canonical modules}

The proof of Theorem~\ref{t:reciprocal} relies on the interpretation
\begin{eqnarray*}
  F_{C\minus\Delta}(x) &=& \Hilb( I_\Delta, x)
\end{eqnarray*}
of the lattice point enumerator as the multigraded Hilbert series
$\Hilb(M,x)$ of the $\ZZ^d$-graded module~$I_\Delta$.  The proof
emphasizes the relations between between $\kk[\Delta], I_\Delta,$
and~$I_{\Delta'}$ by taking homomorphisms into the canonical module.
Throughout we will freely use concepts from combinatorial commutative
algebra that may be found in \cite[Chapter~6]{BrunsHerzog},
\cite[Chapter~1]{Stanley-greenbook}, or
\cite[Part~II]{MillerSturmfels}.

Hochster~\cite{Hochster} showed that the semigroup ring~$\kk[Q]$ is
Cohen--Macaulay whenever $Q$ is saturated.  For a graded
Cohen--Macaulay ring $R$ of dimension $d$, there is the notion of its
{\em canonical module}\/~$\omega_R$.  For $R=\kk[Q]$, it is known
(see e.g. \cite[\S I.13]{Stanley-greenbook}, \cite[\S 13.5]{MillerSturmfels}) that
the canonical module~$\omega_{\kk[\Delta]}$ is the ideal in $\kk[Q]$
spanned $\kk$-linearly by the monomials whose exponents lie in the
interior of the cone~$C$.  Given a Cohen--Macaulay ring~$R$ of
dimension~$d$, and $M$ a Cohen--Macaulay $R$-module of dimension~$e$,
one can \mbox{define the {\em canonical module}\/ of~$M$~by}
\begin{eqnarray*}
  \omega_R(M) &:=& \Ext^{d-e}_R(M, \omega_R).
\end{eqnarray*}
Graded local duality implies that $\omega_R(M)$ is again a
Cohen--Macaulay $R$-module of dimension~$e$, and that
$\omega_R(\omega_R(M)) \cong M$ as $R$-modules.

\begin{proposition} \label{CM-and-canonical-modules}
Let $\Delta \subset C$ be a subcomplex of dimension~$d$, and set $Q =
C \cap \ZZ^d$.
\begin{numbered}
\item
$\kk[\Delta]$ is Cohen--Macaulay if and only if $I_\Delta$ is a
Cohen--Macaulay\/ $\kk[Q]$-module.
\item
$I_\Delta$ is Cohen--Macaulay if and only if $I_{\Delta'}$ is
Cohen--Macaulay, and in this case there is an isomorphism $I_{\Delta'}
\cong \omega_{\kk[Q]}(I_\Delta)$ as $\kk[Q]$-modules.
\end{numbered}
\end{proposition}

\begin{proof}
For the first assertion we use the fact that a graded module~$M$
over~$\kk[Q]$ is Cohen--Macaulay if and only if its local cohomology
$H^i_\mm(M)$ with respect to the maximal ideal $\mm=\<x^a \mid 0 \neq
a \in Q\>$ vanishes for $i$ in the range $[0,\dim(M)-1]$.

The short exact sequence $0 \rightarrow I_\Delta \rightarrow \kk[Q]
\rightarrow \kk[\Delta] \rightarrow 0$ gives a long exact local
cohomology sequence containing the four term sequence
\begin{equation} \label{long-exact-snippet}
H^i_\mm(\kk[Q])
  \rightarrow H^i_\mm(\kk[\Delta]) 
  \rightarrow H^{i+1}_\mm(I_\Delta)
  \rightarrow H^{i+1}_\mm(\kk[Q]).
\end{equation}
Cohen--Macaulayness of $\kk[Q]$ implies that the two outermost terms
of \eqref{long-exact-snippet} vanish for $i$ in the range $[0,d-2]$,
so that the middle map is an isomorphism.  As $\kk[\Delta]$ has
dimension $d-1$, it is Cohen--Macaulay if and only if
$H^i_\mm(\kk[\Delta])$ vanishes for $i \in [0,d-2]$.  As $I_\Delta$
has dimension $d$, it is Cohen--Macaulay if and only if
$H^{i+1}_\mm(I_\Delta)$ vanishes for $i \in [-1,d-2]$.  Noting that
$H^0_\mm(I_\Delta)$ always vanishes due to the fact that $I_\Delta$ is
torsion-free as a $\kk[Q]$-module, the first assertion follows.

For the second assertion, assuming that $I_\Delta$ is Cohen--Macaulay,
we prove a string of easy isomorphisms and equalities:
\begin{equation} \label{equality-string}
\begin{array}{rcl}
I_{\Delta'}  &=& \left( \omega_{\kk[Q]}:I_\Delta \right) \\
             &\cong& \Hom_{\kk[Q]}(I_\Delta, \omega_{\kk[Q]}) \\
             &=& \Ext^0_{\kk[Q]}(I_\Delta, \omega_{\kk[Q]}) \\
             &=& \omega_{\kk[Q]}(I_\Delta)
\end{array}
\end{equation}
in which $(J:I)=\{r \in R: rI \subset J\}$ is the {\em colon ideal}\/
for two ideals $I,J$ in a ring $R$.

The last two equalities in \eqref{equality-string} are essentially
definitions.  To prove the first equality, we claim that if $x^a \in
I_\Delta'$ and $x^b \in I_\Delta$, then $x^a \cdot x^b = x^{a+b} \in
\omega_{\kk[Q]}$.  Using the linear inequalities from
Section~\ref{main}, this holds because
\begin{bullets}
\item
$\ell_F(a) \geq 0$ and $\ell_F(b) \geq 0$ for all facets $F \in \FF$,
\item
$\ell_F(a) > 0$ for $F \in \GG$,
\item
$\ell_F(b) > 0$ for $F \in \FF-\GG$.
\end{bullets}
Thus $I_{\Delta'} \subset \left( \omega_{\kk[Q]}:I_\Delta \right)$.  The
reverse inclusion follows by a similar argument.

The isomorphism in the second line of \eqref{equality-string} 
follows from a general fact:  for any
two $\ZZ^d$-graded ideals $I, J$ in $\kk[Q]$, one has
\begin{eqnarray*}
  \Hom_{\kk[Q]}(I,J) &\cong& (J:I).
\end{eqnarray*}
To prove this, assume $\phi:I \rightarrow J$ is a $\kk[Q]$-module
homomorphism that is $\ZZ^d$-homogeneous of degree~$c$.  Since each
$\ZZ^d$-graded component of $I$ or $J$ is a $\kk$-vector space of
dimension at most $1$, for each monomial~$x^a$ in~$I$ there exists a
scalar $\lambda_a \in \kk$ such that \mbox{$\phi(x^a)=\lambda_a
x^{a+c}$}.  We claim that these scalars $\lambda_a$ are all equal to a
single scalar $\lambda$.  Indeed, given $x^a, x^b$ in~$I$, the fact
that $\phi$ is a $\kk[Q]$-module homomorphism forces both $\lambda_a =
\lambda_{a+b}$ and $\lambda_b=\lambda_{a+b}$.  Thus for some $\lambda
\in \kk$, one has $\phi=\lambda \cdot \phi_c$, where
$\phi_c(x^a)=x^{a+c}$.  Furthermore, if $\lambda \neq 0$ then $x^c \in
(J:I)$.  We conclude that the map
\begin{eqnarray*}
  (J:I) &\rightarrow& \Hom_{\kk[Q]}(I,J)\\
   x^c  &\mapsto& \phi_c
\end{eqnarray*}
is an isomorphism of $\kk[Q]$-modules.
\end{proof}

\begin{trivlist}
\item {\it Proof of Theorem~\ref{t:reciprocal}.\,}
The key fact (see \cite[\S I.12]{Stanley-greenbook}, for instance) is
that for any Cohen--Macaulay $\kk[Q]$-module~$M$ of dimension~$d$,
\begin{eqnarray*}
  \Hilb(\omega_{\kk[Q]}(M);x^{-1}) &=& (-1)^d \Hilb(M;x).
\end{eqnarray*}
Therefore when $I_\Delta$ is Cohen--Macaulay,
Proposition~\ref{CM-and-canonical-modules} gives
\begin{eqnarray*}
F(C\minus\Delta'; x^{-1}) &=& \Hilb(I_{\Delta'}; x^{-1})\\
                     &=& (-1)^d \Hilb(I_\Delta; x)\\
                     &=& (-1)^d F(C\minus\Delta; x).
\makebox[0pt][l]{\makebox[.265\linewidth][r]{\qed}}
\end{eqnarray*}
\end{trivlist}

\section{Cohen--Macaulay $d$-complexes in $\RR^d$} \label{d-space}

The goal of this section is to prove
Theorem~\ref{CM-d-complexes-in-d-space}.  In this section, $K$ will be
a finite polyhedral complex embedded piecewise linearly in $\RR^d$.
That is, $K=\{P_i\} = \{P_i\}_{i\in I}$ is a finite collection of
convex polytopes in~$\RR^d$ containing the faces of any polytope~$P_i$
in~$K$, and for which any two polytopes $P_i$ and $P_j$ in~$K$
intersect in a common (possibly empty) face of each.

It will be convenient to pass between \PL-embeddings of such
polyhedral complexes into $\RR^d$, and \PL-embeddings into the
boundary of a $(d+1)$-polytope.  In one direction, this passage is
easy, as we now show.

\begin{proposition} \label{p:Schlegel}
Let $K$ be a finite polyhedral $d$-dimensional complex \PL-embedded as
a proper subset of the boundary of a $(d+1)$-polytope $P$ (but not
necessarily as a subcomplex of the boundary).  Then $K$ has a
\PL-embedding into~$\RR^d$.
\end{proposition}
\begin{proof}
We first reduce to the case where $K$ avoids at least one facet of~$P$
entirely.  Since $K$ is a compact proper subset of the boundary of
$P$, there exists at least one facet~$F$ of~$P$ whose interior is not
contained in~$K$.  Let~$\sigma$ be a $d$-dimensional simplex
\PL-embedded in the complement $F \minus K$, and let $P'$ be a
$(d+1)$-simplex obtained by taking the pyramid over~$\sigma$ whose
apex is any interior point of~$P$.  Then projecting~$K$ from any
interior point of~$P'$ onto the boundary of~$P'$ gives a \PL-embedding
of~$K$ into this boundary, avoiding the facet~$\sigma$ of~$P'$
entirely.

Once $K$ avoids a facet $F$ of $P$ entirely, it is \PL-homeomorphic to
a subcomplex of a {\em Schlegel diagram}\/ for~$P$ in~$\RR^d$
\cite[Definition~5.5]{Ziegler} with $F$ as the bounding facet.%
\end{proof}

For the other direction, we use a construction of
J.\thinspace{}Shewchuk.

\begin{theorem}[Shewchuk] \label{t:shewchuk}
Let $K$ be a polyhedral complex \PL-embedded in~$\RR^d$.  Then $K$ is
\PL-homemorphic to a subcomplex of the boundary of a $(d+1)$-polytope.
\end{theorem}
\begin{proof}
Consider an arrangement $\AAA=\{H_i\}$ of finitely many affine
hyperplanes in~$\RR^d$ with the property that every polytope~$P$
in~$K$ is an intersection of closed halfspaces bounded by some subset
of the hyperplanes in~$\AAA$; since~$K$ contains only finitely many
polytopes, such arrangements exist.

Let $K'$ be the subdivision of~$K$ induced by its intersection with
the hyperplanes of~$\AAA$, so that $K'$ is a finite subcomplex of the
polyhedral subdivision~$\hat{K}$ of~$\RR^d$ induced by~$\AAA$.  Then
$\hat{K}$ is a {\em regular}\/ (or {\em coherent}\/
\cite[Definition~7.2.3]{GKZ}) subdivision; its cells are exactly the
domains of linearity for the piecewise-linear convex function
$$
\begin{array}{rcl}
  f:\RR^d &\rightarrow & \RR \\
      x   & \mapsto    & \sum_{i} d(x, H_i)
\end{array}
$$
in which $d(x,H)$ denotes the (piecewise-linear, convex) function
defined as the distance from $x$ to the affine hyperplane $H$.  Since
$K$ is finite, $f$ achieves a maximum value, say $M$, on $K$. Then for
any $\epsilon > 0$, the $(d+1)$-dimensional convex polytope
$$
  \{(x,x_{d+1}) \in \RR^{d+1} \mid f(x) \leq x_{d+1} \leq M+\epsilon\}
$$
contains the graph
$$
  \{(x,f(x)) \mid x \in K\}
$$
of the restricted function $f|_{K}$ as a polyhedral subcomplex of its
lower hull.  Furthermore, the projection $\RR^{d+1} \rightarrow \RR^d$
gives an isomorphism of this subcomplex onto~$K'$.%
\end{proof}

Theorem~\ref{CM-d-complexes-in-d-space}.2 will
follow immediately from the following construction of B. Mazur, which
is famous in the topology community.

\begin{proposition} \label{Mazur-construction}
$\RR^4$ contains a \PL-embedded finite simplicial complex~$K$ that
triangulates a contractible $4$-manifold, but whose boundary is not
simply-connected.  In particular, $K$ is Cohen--Macaulay over every
field\/~$\kk$ but not \mbox{homeomorphic to a $4$-ball}.
\end{proposition}
\begin{proof}
Mazur \cite[Corollary 1]{Mazur} constructs a finite simplicial
complex~$K$ that is contractible, has non-simply-connected
boundary~$\partial K$, and enjoys the further property that its
``double''~$2K$ (obtained by identifying two disjoint copies of~$K$
along their boundaries) is \PL-isomorphic to the boundary of a $5$-cube.  Thus 
$K$ is \PL-embedded as a proper subset of the boundary of a $5$-polytope, and
hence has a \PL-embedding in $\RR^4$ by Proposition~\ref{p:Schlegel}.
\end{proof}

\begin{proofof}{Theorem~\ref{CM-d-complexes-in-d-space}.2.}
This is a consequence of Theorem~\ref{t:shewchuk} and
Proposition~\ref{Mazur-construction}.%
\end{proofof}

We next turn to Theorem~\ref{CM-d-complexes-in-d-space}.1.  
Fix a field $\kk$.  For each
nonnegative integer~$d$, consider two related assertions $A_d$
and~$A'_d$ concerning finite $d$-dimensional polyhedral complexes,
where we write `CM' for `Cohen--Macaulay over $\kk$'.
\begin{enumerate}
\item[$A_d$:] Every \PL-embedded CM $d$-complex in~$\RR^d$ is
homeomorphic \mbox{to a $d$-ball}.
\item[$A'_d$:] Every \PL-embedded CM $d$-complex in~$\RR^d$ is a
$d$-manifold with boundary.
\end{enumerate}

Assertion $A_d$ is false for $d=4$, as shown by
Proposition~\ref{Mazur-construction}.  We wish to show it is true for
$d \leq 3$, as this would in particular prove the assertion of
Theorem~\ref{CM-d-complexes-in-d-space}.1.

Throughout the remainder of this section, all homology and cohomology
groups are reduced, and taken with coefficients in~$\kk$.  We will
also use implicitly without further mention the fact that any
Cohen--Macaulay $d$-complex $K$ embedded in $\RR^d$ must necessarily
be $\kk$-acyclic: the Cohen--Macaulay hypothesis gives $H_i(K)=0$ for
$i<d$, and Alexander duality within the one-point compactification of
$\RR^d$ implies~\mbox{$H_d(K)=0$}.

In the proof of the next lemma, we use the notion of {\em links}\/
(sometimes also called {\em vertex figures}\/) of faces
(polytopes)~$F$ in a polyhedral complex~$K$ that is \PL-embedded
in~$\RR^d$.  For each face~$F$, we (noncanonically) construct a
polyhedral complex $\link_K(F)$ that models the link.  First, write
$\RR^d/F$ for the quotient of~$\RR^d$ by the unique linear subspace
parallel to the affine span of~$F$.  Then choose a small simplex
$\sigma$ containing the point $F/F \in \RR^d/F$ in its interior.  Each
face~$G$ of~$K$ containing~$F$ has an image~$G/F$ in~$\RR^d/F$ whose
intersection with each face of~$\sigma$ is a polytope.  These
polytopes constitute the faces of a polyhedral complex \PL-embedded in
the boundary of~$\sigma$, and we take $\link_K(F)$ to be this complex.

\begin{lemma} \label{manifold-lemma}
If assertion $A_\delta$ holds for every $\delta < d$ then
assertion $A'_d$ holds.
\end{lemma}
\begin{proof}
Assume that $K$ satisfies the hypotheses of $A'_d$.  To show that $K$
is a $d$-manifold with boundary, it suffices to show that $\link_K(F)$
is either a $\delta$-sphere or a $\delta$-ball for every
$e$-dimensional face~$F$ of~$K$, where $\delta = d-e-1$.  We use the
fact that the link of any face in any Cohen--Macaulay $d$-complex is
Cohen--Macaulay.  This holds in our case because Cohen--Macaulayness
is a topological property (see Section~\ref{main}), so we can
barycentrically subdivide and use the corresponding fact for
simplicial complexes (which follows from Reisner's criterion for
Cohen--Macaulayness via links \cite[\S II.4]{Stanley-greenbook}).

By construction, $\link_K(F)$ is a $\delta$-dimensional polyhedral
complex \PL-embedded in the boundary of 
the small $(\delta+1)$-simplex $\sigma$ around $F/F$.  If the barycenter of
$F$ is an interior point of the manifold $K$, then $\link_K(F)$
is a polyhedral subdivision of the entire (topologically $\delta$-spherical) 
boundary of $\sigma$;  otherwise it is embedded as a proper subset.
In the latter case, Proposition~\ref{p:Schlegel} and assertion~$A_\delta$ 
apply to show that $\link_K(F)$ is a $\delta$-ball, as desired.%
\end{proof}

\begin{theorem} \label{A0-A1-A2-A3}
Assertion $A_d$ holds for $d \leq 3$.
\end{theorem}
\begin{proof}
Assertions $A_0, A_1$ are trivial.  Together they imply
assertion~$A_2'$ via Lemma~\ref{manifold-lemma}.  {}From this,
deducing the stronger assertion $A_2$ is a straightforward exercise
using
\begin{bullets}
\item
the fact that the boundary $ \partial K$ is a disjoint union of
$1$-spheres (possibly nested) embedded in $\RR^2$,
\item
the Jordan Curve Theorem, and 
\item
$H_0(K)=H_1(K)=0$.
\end{bullets}

To prove~$A_3$, we may assume~$A'_3$ by Lemma~\ref{manifold-lemma},
and hence assume that $K$ is a Cohen--Macaulay $3$-manifold with
boundary, embedded in~$\RR^3$.  Thus~$H_1(K)=0$, and hence
Lemma~\ref{li-lemma} below forces~\mbox{$H_1(\partial K)=0$}.  Since
$\partial K$ is orientable, this implies that $\partial K$ is a
disjoint union of (possibly nested) $2$-spheres.  It is then another
straightforward exercise using the Jordan--Brouwer Separation Theorem,
along with the fact that $H_0(K)=H_2(K)=0$, to deduce that $\partial
K$ must consist of a single $2$-sphere, with $K$ its interior.  The
Alexander--Schoenflies Theorem then implies that $K$ is a $3$-ball.%
\end{proof}

\begin{proofof}{Theorem~\ref{CM-d-complexes-in-d-space}.1.}
Immediate from Theorem~\ref{A0-A1-A2-A3} and
Proposition~\ref{p:Schlegel}.
\end{proofof}

The authors thank T.-J. Li for pointing out the following lemma
and proof (cf. \cite[proof of Theorem 6.40]{Vick}), which was used in
the proof of Theorem~\ref{A0-A1-A2-A3}.

\begin{lemma} \label{li-lemma}
For any compact $3$-manifold~$K$ with boundary~$\partial K$,
\begin{eqnarray*}
  \dim_\kk H_1(K;\kk) &\geq& \frac{1}{2} \dim_\kk H_1(\partial K;\kk).
\end{eqnarray*}
\end{lemma}
\begin{proof}
Consider the following diagram, in which the two squares commute:
\begin{equation*}
\begin{CD}
  \Hom(H_1(K),\kk)   @>{\Hom(i_*,\kk)}>>  \Hom(H_1(\partial K),\kk)\\
    @A{}AA                              @AA{}A     \\
       H^1(K)        @>{i^*}>>           H^1(\partial K)\\
   @V{}VV                              @VV{}V     \\
    H_2(K,\partial K)@>{j_*}>>   H_1(\partial K)  @>{i_*}>> H_1(K)   &\\
\end{CD}
\end{equation*}
The vertical maps are all isomorphisms.  The two vertical maps in the
top square come from the universal coefficient theorem relating
cohomology and homology with coefficients in~$\kk$.  The two vertical
maps in the bottom square are duality isomorphisms, the left coming
from Poincar\'e--Lefschetz duality for $(K, \partial K)$ and the right
from Poincar\'e duality for~$\partial K$.

The inclusion $\partial K \overset i\hookrightarrow K$ induces three
of the horizontal maps.  The last row is exact at its middle term,
forming part of the long exact sequence for the pair $(K,\partial K)$,
in which $j_*$ is a connecting homorphism.  Thus
$$
\begin{array}{rcccccccl}
  \nullity(i_*) &=& \rank(j_*) &=& \rank(i^*) &=& \rank(\Hom(i_*,\kk))
  &=& \rank(i_*).
\end{array}
$$
On the other hand, 
$$
\dim_\kk H_1(\partial K) \,\, =\,\, \rank(i_*) + \nullity(i_*) 
                         \,\,= \,\,2 \, \rank(i_*) 
                         \,\,\leq \,\, 2 \dim_\kk H_1(K),
$$
which completes the proof.
\end{proof}

\section*{Acknowledgements}
The authors thank Jonathan Shewchuk for allowing his construction in
Theorem~\ref{t:shewchuk} to be included here, and to Robion Kirby for
pointing out Mazur's construction.  Tian-Jun Li kindly provided us
with crucial technical help, and Lou Billera, Anders Bj\"orner, Don
Kahn, and Francisco Santos gave helpful comments.

\end{document}